\title{ ~~\\ Primitive root producing quadratics}
\author{Pieter Moree}
\documentclass[12pt]{article}
\usepackage{amssymb, latexsym, amsfonts}
\def\@ptsize{2}
\newtheorem{Thm}{Theorem}
\newtheorem{Def}{Definition}
\newtheorem{Heuristic}{Heuristic}
\newtheorem{Lem}{Lemma}
\newtheorem{Problem}{Problem}
\newtheorem{Conjecture}{Conjecture}
\newtheorem{Cor}{Corollary}

\newtheorem{Prop}{Proposition}
\newcommand{\qed}{\hfill $\Box$}

\begin{document}
\date{}
\maketitle {\def\thefootnote{} \footnote{\noindent P. Moree:
Max-Planck-Institut f\"ur Mathematik, Vivatsgasse 7, D-53111 Bonn,
Deutschland, e-mail: moree@mpim-bonn.mpg.de}} {\def\thefootnote{}
}
\begin{abstract}
\noindent D.H. Lehmer found a quadratic polynomial such that $326$
is a primitive root for the first $206$ primes represented by
this polynomial. It is shown that this is related to the class
number one problem and prime producing quadratics. More
impressive examples in the same spirit are given using recent
results on prime producing quadratics. Y. Gallot holds the current
record in which $206$ is being replaced by $31082$.
\end{abstract}
\section{Introduction}
In their celebrated book Ireland and Rosen \cite{IR} write (p.
47): `Lehmer discovered the following curious result. The first
prime of the form $326n^2+3$ for which $326$ is not a primitive
root must be bigger than $10$ million. He mentions other results
of the same nature. It would be interesting to see what is
responsible for this strange behavior'. Using e.g. Maple one
easily checks that $326$ is a primitive root mod $p$ for the first
$206$ primes of the form $326n^2+3$ (they satisfy $0\le n\le
2374$), but is not for $p=1838843753=326\cdot
2375^{2}+3$.\\
\indent Note that $326=2\cdot 163$ and recall that the class
number of $\mathbb Q(\sqrt{-163})$ equals one. It will be shown in
this note that there is a connection between this fact, finding
prime producing polynomials and Lehmer's observation. This
suggests the (apparently) unexplored idea of finding `primitive
root producing polynomials'. We say a polynomial $f(X)$ is {\it
primitive root producing} if for a prescribed integer $g$, $g$ is
very frequently a primitive root modulo those primes that are
assumed as values by $f$. We are especially interested in
quadratics $f$ such that a given integer $g$ is a primitive root
modulo $p$ for as many consecutive $n$ as possible for which
$f(n)$ is prime.
\begin{Def} Given integers $g$ and $f(X)\in \mathbb Z[X]$, let
$p_1(g,f),p_2(g,f),\dots$ be the consecutive primes of the form
$f(n)$ with $n\ge 0$ that do not divide $g$. We let $r$ be
the largest integer $r$ (if this exists) such that $g$ is a
primitive root mod $p$ for all primes $p_j(g,f)$ with $1\le j\le
r$. We let $c_g(f)$ be the number of distinct primes amongst
$p_j(g,f)$ with $1\le j\le r$.
\end{Def}
Thus, for example, $c_{326}(f)=206$, with $f(X)=326X^2+3$.
\begin{Problem}
\label{P1} Find $g$ and $f$ such that $c_g(f)$ is as large as
possible.
\end{Problem}
By the Chinese Remainder Theorem we know, that given any finite
set of odd primes one can find $g$ such that $g$ is a primitive
root modulo for each of these primes, thus one should require $g$
to be small in comparison with the coefficients of $f$. We say $g$
is small in this context if $|g|<10^{c_g(f)/3}$ (see
Section \ref{likeli} for an explanation).\\
\indent The starting point of Lehmer's paper was a letter he
received in 1957 from one Raymond Griffin (then living in Dallas,
Texas). In this letter Griffin suggested that the decimal
expansions of $1/p$ should have period length $p-1$ for all primes
of the form $10n^2+7$. Note that, if $p\nmid 10$, ord$_p(10)=v$
iff the period of the decimal expansion of $1/p$ is $v$. The first
$16$ primes $p$ of the form $10n^2+7$ have indeed decimal period
$p-1$, but this is not true for $p=7297$, the $17$th such prime.
Griffin's conjecture suggests the following problem:
\begin{Problem}
Given a prescribed integer $g$ in
$$G:=\{g\in \mathbb Z:g\ne -1{\rm ~and~}g\ne b^2,~b\in \mathbb Z\},$$ find a
quadratic polynomial $f$ such that $c_g(f)$ is as large as
possible.
\end{Problem}
Note that an integer in $\mathbb Z\backslash G$ is a primitive
root for only finitely many primes. Since Problem 2 is an easy variant
of Problem 1, we will not discuss this further.\\
\indent Of course there is no need to
restrict to quadratic polynomials, but this is what we shall
do in this paper. Since at present it is not even known whether $n^2+1$ is
prime infinitely often, we can only expect to gain some insight on
assuming certain conjectures. In the next section we
briefly recall some relevant conjectures.

\section{Prerequisites on two conjectures}
Let $f(X)$ be an irreducible polynomial of content $1$ in $\mathbb
Q[X]$  with integer coefficients. By a special case of a conjecture
due to Bateman and Horn
\cite{BH} $\pi_f(x)$, the
number of integers $0\le n\le x$ such that $f(n)$ is prime, should
satisfy, as $x$ tends to infinity,
$$\pi_f(x)\sim {H(f)\over {\rm deg}(f)}{x\over \log x},{\rm ~where~}
H(f)=\prod_p{1-{N_p(f)\over p}\over 1-{1\over p}},$$ and $N_p(f) =
\#\{n({\rm mod~}p):f(n)\equiv 0({\rm mod~}p)\}$. We say a
congruence class modulo an integer $m$ is {\it allowable} if for
any number $r$ in it we have $(f(r),m)=1$ and thus, e.g.,
$p-N_p(f)$ denotes the number of
allowable congruence classes modulo $p$.\\
\indent If we fix the degree of $f$ then, by the fundamental lemma
of the sieve, we have uniformly in $f$, that
\begin{equation}
\label{tweede}
\pi_f(x)\ll \prod_{p\le x}(1-{N_p(f)\over p})x.
\end{equation}
If $\sum_{p>x}(1-N_p(f))/p\ll 1$, then (\ref{tweede}) becomes
$\pi_f(x)\ll H(f)x/{\rm deg}(f)\log x$, uniformly \cite{GM}.\\
\indent Let $\cal F$ be the set of quadratic polynomials
$aX^2+bX+c$ with $a>0$, $b,c$ integers such that gcd$(a,b,c)=1$,
$d=b^2-4ac$ is not a square and $a+b$ and $c$ are not both even.
Then, as $x$ tends to infinity, Hardy-Littlewood's Conjecture F \cite{HL},
a special case of the Bateman-Horn conjecture, asserts that
\begin{equation}
\label{HLL} \pi_f(x)\sim \epsilon {x\over \log x}\prod_{p>2\atop
p|(a,b)}{p\over p-1} \prod_{p>2\atop p\nmid a}\left(1-{({d\over
p})\over p-1}\right),
\end{equation}
where $\epsilon=1$ if $a+b$ is even and $\epsilon=1/2$ otherwise. 
For $f\in {\cal F}$ it is easily shown that 
\begin{equation}
\label{afschatting1}
{a\over \varphi(a)L(1,(d/.)}\ll H(f)\ll {a\over \varphi(a)L(1,(d/.)}.
\end{equation}
\indent For our purposes the following weaker conjecture, which is
implied by Hardy-Littlewood's Conjecture F, will suffice.
\begin{Conjecture}
\label{equidi} Let $m\ge 2$ be an integer. Suppose that $f(X)\in
\mathbb Z[X]$ represents infinitely many primes, then the $n$ for
which $f(n)$ is prime are asymptotically equidistributed over the
allowable congruence classes modulo $m$.
\end{Conjecture}
\indent Finally we recall the prime $k$-tuplets conjecture
(TC($k$)). This conjecture seems to be due to Dickson (1904).
\begin{Conjecture}
Let $k\ge 1$ and let $A_1,\dots,A_k,B_1,\dots,B_k$ be integers
with $A_j>0$ for $j=1,\dots,k$. Suppose that for each prime $p$
there exists an integer $n_p$ such that $p$ does not divide
$\prod_{j=1}^k (A_jn_p+B_j)$, then there exist infinitely many
integers $n$ such that $A_jn+B_j$ is prime for $1\le j\le k$.
\end{Conjecture}

\section{On the likelihood of finding $c_g(f)=m$}
\label{likeli}

\noindent Given a finite set of primes $\{p_1,\dots,p_s\}$ let
$P=\prod_{i=1}^s p_i$. There are $\prod_{i}^s \varphi(p_i-1)$
residue classes modulo $P$ such that if $g$ is in any of them it
is a primitive root for every prime dividing $P$. Assuming
equidistribution we expect that the smallest of them is roughly of
size $Q:=\prod_{i=1}^s (p_i-1)/\varphi(p_i-1)$. It is an easy
exercise in analytic number theory to evaluate the average value
of $(p-1)/\varphi(p-1)$. To this end note that
$$\sum_{p\le x}{p-1\over \varphi(p-1)}=\sum_{d\le x}{\mu(d)^2\over
\varphi(d)}\pi(x;d,1),$$ where $\pi(x;d,1)$ denotes the number of
primes $q\le x$ such that $q\equiv 1({\rm mod~}d)$. Proceeding as
in the proof of Lemma 1 of \cite{MoreeII} one then finds that for
every $C>1$ one has
$$\sum_{p\le x}{p-1\over \varphi(p-1)}=
B{\rm Li}(x)+O({x\over \log^{C}x}),~{\rm with~}B=\prod_{q{\rm
~prime}}\left(1+{1\over (q-1)^2}\right),$$ where the implied
constant may depend on $C$ and Li$(x)$ denotes the logarithmic
integral. This improves on an estimate due to Murata \cite{gp}.
Expressing $B$ in terms of zeta values, cf. \cite{C}, one finds
$B=2.826419997067\dots$ . Thus $Q$ is roughly of size
$B^s\approx 10^{0.45s}$. This motivates the definition of small $g$ in
the introduction.\\
\indent Likewise one can wonder about the probability that a given
$g$ is a primitive root for our finite set of primes. An estimate
for this is given by $1/Q$ and should be roughly $10^{-0.45s}$.
Thus a measure for the likelihood of having $c_g(f)=m$ (by random choice
of $f$
and $g$) is $10^{-m/2}$.

\section{Lehmer's observation}
The following trivial result will play an important role in the
explanation of Lehmer's observation (and in finding some more
impressive variants of it):
\begin{Lem}
\label{een} Let $\alpha\ge 0$ be an integer. Let $p$ be a prime
and $g$ an integer coprime with $p$. Define $r_{p}(g):=[(\mathbb
Z/p\mathbb Z)^*:\langle g\rangle]$ (the residual index of $g({\rm
mod~}p)$). Let $d_1,d_2$ be positive integers. Let $p$ be a prime
of the form $2^{\alpha}d_1n^2+d_22^{\alpha}+1$. If $q$ is an odd
prime with $({-d_1d_2\over q})\ne 1$ and $q\nmid d_2$, then
$q\nmid r_p(g)$.
\end{Lem}
{\it Proof}. The equation $2^{\alpha}d_1X^2+d_22^{\alpha}+1=1$ is
solvable mod $q$ if and only if $({-d_1d_2\over q})=1$ or $q|d_2$.
Since by assumption $({-d_1d_2\over q})\ne 1$ and $q\nmid d_2$, it
follows that $p\not\equiv 1({\rm mod~}q)$. From this and
$r_p(g)|p-1$, it then
follows that $q\nmid r_p(g)$. \qed\\

\noindent Using Lemma \ref{een} it is easy to deduce the following
proposition.
\begin{Prop}
\label{prop} Let $k$ be a non-zero integer. Let $g\in
\{-163,-3,6,326\}$. If $p$ is a prime not dividing $kg$ and
$p=326n^2+3$, then $(r_p(k^2g),2\cdot 3\cdots 37)=1$.
\end{Prop}
{\it Proof}. Using quadratic reciprocity one deduces that
$({k^2g\over p})=-1$ and hence $2\nmid r_p(k^2 326)$. Let $q$ be
an odd prime not exceeding $37$. It is easy to check (using e.g.
quadratic reciprocity) that $({-163\over q})=-1$ and thus, by
Lemma \ref{een}, $q\nmid
r_p(k^2g)$. \qed\\

\noindent Put $L(X)=326X^2+3$. The latter result shows that if
$326$ is not a primitive root modulo a prime $p=L(n)$, then
$r_p(326)\ge 41$. Since this is rather unlikely to happen, we
expect to find a reasonably long string of primes of the form
$L(n)$ before we find a prime $p$ for which $326$ is not a
primitive root mod $p$. This is precisely what happens: we have to
wait until $n=2375$ and hence $p=1838843753$, for $326$ not
to be a primitive root mod $p$ (we have $r_p(326)=83$).\\
\indent Supposing $p=L(n)$ to be prime, one can wonder about the
probability that $r_{p}(326)>1$. For this to happen $r_{p}(326)$
must be divisible by some odd prime $q$ such that $({-163\over
q})=1$. In this case $n$ has to be in one of two residue classes
mod $q$ and, moreover, we need to have $326^{p-1\over q}\equiv
1({\rm mod~}p)$. Since $326^{p-1\over q}$ is merely one out of the
$q$ solutions of $x^q\equiv 1({\rm mod~}p)$, one heuristically
expects that $326^{p-1\over q}\equiv 1({\rm mod~}p)$ with
probability $1/q$. We thus expect that with probability
\begin{equation}
\label{lehmertje}
\prod_{({-163\over q})=1}\left(1-{2\over q^2}\right)=0.99337\ldots
\end{equation}
a prime of the form $p=L(n)$ will have $326$ as a primitive
root. This argument is taken from Lehmer's paper. He implicitly
assumes that the $n$ for which $f(n)$ is prime are asymptotically equally
distributed over the congruence classes modulo $q$, instead of
over the allowable congruence classes modulo $q$. On correcting
for this one arrives at a probability of
\begin{equation}
\label{p1}
p_1:=\prod_{({-163\over q})=1}\left(1-{2\over q(q-1-({-978\over
q}))}\right)=0.99323\ldots.
\end{equation}
For $0\le n\le 5\cdot 10^6$ there
are $240862$ primes $p=L(n)$ of which $239239$ have $326$ as a
primitive root. Note that $239239/240862\approx 0.99326\ldots $ .\\
\indent Instead of taking $326$ as base, Proposition \ref{prop}
suggests we could take $k^2 326$ as a base and vary over $k$. Assuming that 
each prime $p=L(n)$ has a probability $p_1$
of having $k^2 326$ as a primitive root we might expect that
$$\lim_{x\rightarrow \infty}{1\over x}\sum_{k\le x}g_{k^2326}(f)\approx
\sum_{j=1}^{\infty}jp_1^j(1-p_1)={p_1\over 1-p_1},$$ that is
equals about $150$ (note that the `probability' that
$g_{k^2326}(f)=j$ equals $p_1^j-p_1^{j+1}=p_1^j(1-p_1)$). For
$k\le 5000$ it turns out that the average is around $180$. Note
that in the averaging process there is a very strong bias
towards the smallest primes of the form $p=L(n)$. This might explain the observed discrepancy.\\
\indent The most interesting quantity for our purposes is
$\max_{1\le k\le s}g_{k^2 326}(L)$. For this one expects the
outcome
$$M(p_1,s):=\sum_{j=1}^{\infty}j\left((1-p_1^{j+1})^s-
(1-p_1^{j})^s\right).$$ It is not difficult to show that, as $s$
tends to infinity,
\begin{equation}
\label{oudprobleempje} M(p_1,s)\sim {\log s\over \log(1/p_1)},
\end{equation}
and that this holds more generally for any value of $p_1$
satisfying $0<p_1<1$ \cite{Moree}. By more subtle techniques
\cite{BvV, vV} it can be shown that
$$M(p_1,s)\approx {1\over \log(1/p_1)}\sum_{r=1}^s{1\over r}-{1\over 2},$$
where the approximation is remarkably good and $0<p_1<1$. The
interpretation of the latter result is somewhat disappointing: if
one has found $M(s):=\max_{1\le k\le s}g_{k^2 326}(L)$ with
$s=10^6$, say, then in order to find a $k$ such that $g_{k^2
326}(L)\ge 2M$ one expects to have to compute $g_{k^2 326}(L)$ for
all $k$ up to around $10^{12}$ in order to achieve this. The
numerics seem to confirm the slow growth of $M(s)$. For example,
$M(350)=1123$ and
$M(25000)=1614$.\\
\indent One can wonder how `special' it is to find a given value of
$c_{k^2g}(L)$. An obvious measure for this is the smallest integer
$s$ such that $M(p_1,s)=c_{k^2g}(L)$. For $1614$ for example this
is around $32500$, i.e., one would expect to try around $32500$ values
of $k$ before finding $c_{k^2g}(L)\ge 1614$.\\
\indent Griffin's and Lehmer's polynomial for $g=10$, respectively
$g=326$ show that there are quadratic polynomials $f$ and integers
$g$ such that $({g\over p})\ne 1$ for all primes of the form
$f(n)$, i.e. all the primes $p=f(n)$ are inert in $\mathbb
Q(\sqrt{g})$. In the next section we investigate this situation
further.

\section{On the splitting of primes $p=f(n)$ in a quadratic field}
This section is devoted to a conditional result on the splitting
behaviour of primes of the form $p=f(n)$ in a prescribed quadratic
field $K$. In the case where $f$ is quadratic we will make this
result more
explicit.\\
\indent Let $d>1$ be an odd squarefree integer. Put
\begin{equation}
\label{definitie} a_d(f)={\sum_{r({\rm mod~}d)}\left({f(r)\over
d}\right)\over \#\{r({\rm mod~}d):(f(r),d)=1\}}.
\end{equation}
Note that $-1\le a_d(f)\le 1$. By the Chinese Remainder Theorem
and the multiplicative property of the Jacobi symbol the quantity
$a_d(f)$ is seen to be a multiplicative function on odd squarefree
integers $d$. Thus $a_d(f)=\prod_{p|d}a_p(f)$. Note that if $p>2$
and $N_p(f)$ is even, then $a_p(f)\ne 0$.
\begin{Thm}
\label{thm1} Let $D$ be a fundamental
discriminant. Suppose that $f(n)$ is prime for infinitely many $n$
and that the $n$ for which $f(n)$ is prime are equidistributed
over the residue classes $a({\rm mod ~}D)$ with $(f(a),D)=1$. The
proportion $\tau_D^-(f)$ of primes $p$ satisfying $p=f(n)$ for
some $n$ that are, moreover, inert in a quadratic field of
discriminant $D$ exists and is a rational number. Let $D_1$ be the
largest odd prime divisor of $D$ and assume that $D_1>1$. For
$j=1,3,5$ and $7$ put
$$\alpha_j={\#\{s({\rm mod~}8):f(s)\equiv j({\rm mod~}8)\}\over
4 \#\{s({\rm mod~}2):f(s)\equiv 1({\rm mod~}2)\}}.$$ We have
$$2\tau_D^-(f)=\cases{1-a_{D_1}(f) & if $D$ is odd;\cr
1+(\alpha_3+\alpha_7-\alpha_1-\alpha_5)a_{D_1}(f) & if
$D\equiv 4({\rm mod~}8)$;\cr 
1+(\alpha_3+\alpha_5-\alpha_1-\alpha_7)a_{D_1}(f) & if
$D\equiv 8({\rm mod~}32)$;\cr 
1+(\alpha_5+\alpha_7-\alpha_1-\alpha_3)a_{D_1}(f) & if
$D\equiv 24({\rm mod~}32)$.}$$ Moreover,
$a_{D_1}(f)=\prod_{p|D_1}a_p(f)$, with
$$a_p(f)={\sum_{j=0}^{p-1}({f(j)\over p})\over p-N_p(f)}.$$
\end{Thm}
{\it Proof}. Let us consider the case where $D>1$ and $D\equiv
1({\rm mod~}4)$ first. Note that  $p$ is inert in $K$ iff
$(D/p)=-1$. Since  $D\equiv 1({\rm mod~}4)$, we have $({D\over
p})=({p\over D})$ and thus only the value of $p({\rm mod~}D)$
matters. By assumption the corresponding values of $n$ are
equidistributed asymptotically.  Therefore $\tau_D^-(f)$, the
proportion of primes of the form $f(n)$ which  are inert in $K$,
satisfies
$$\tau_D^-(f)={\#\{r({\rm mod~}D):({f(r)\over D})=-1\}\over \#\{r({\rm
mod~}D):(f(r),D)=1\}}.$$ Let us denote the corresponding proportion of
split primes by $\tau_D^+(f)$. We have $\tau_D^-(f)+\tau_D^+(f)=1$
and $\tau_D^+(f)-\tau_D^-(f)=a_D(f)$, whence
$\tau_D^-(f)=(1-a_D(f))/2$, as required.\\
\indent In case $2|D$ we consider the various congruence classes
modulo $8$ separately. Each of them can then be dealt with as
before (this involves quadratic reciprocity). The details are left
to the
interested reader. \qed\\

\noindent {\tt Remark} 1. Note that under the assumption of
Hardy-Littlewood's Conjecture F the hypothesis of the result is
satisfied. (For then Conjecture \ref{equidi}
holds true.)\\

\noindent {\tt Remark} 2. Notice that the condition that $f(n)$
represents infinitely many primes ensures that $\alpha_j$ exists
for $j=1,3,5$ and $7$. These numbers can be explicitly evaluated,
but this requires a lot of case distinctions.

\subsection{The case where $f$ is quadratic}
Before we state the main result of this section (Proposition
\ref{lemmi}), we need some preliminaries on certain
simple character sums.\\
\indent The following two lemmas are well-known, see \cite[p.
79]{G}. The proof of Lemma \ref{Jacobsthal} given here (suggested
by I. Shparlinski) is more natural than the one in \cite[p.
79]{G}.
\begin{Lem}
\label{Jacobsthal} Let $p$ be an odd prime. Then
$$\sum_{m=0}^{p-1}\left({m^2+a\over p}\right)=\cases{p-1 &if $p|a$;\cr
-1 & otherwise.}$$
\end{Lem}
{\it Proof}. If $p|a$ the assertion is trivial. The result in case
$p\nmid a$ easily follows once we know for how many $y\ne 0$ we
have $m^2+a\equiv y^2({\rm mod~}p)$. Thus we want to have $a\equiv
(y-m)(y+m)({\rm mod~}p)$. Write $u=y-m$ and $v=y+m$. There are
$p-1$ pairs $(u,v)$ satisfying $a\equiv uv({\rm mod~}p)$. Using
that the pairs $(u,v)$ are in bijection with the pairs $(m,y)$,
the proof is then easily completed on distinguishing between the
case $({-a\over p})=-1$ and
$({-a\over p})=1$. \qed\\

\noindent Let $f(x)=ax^2+bx+c$ be a quadratic polynomial. Put
$d=b^2-4ac$ and
$$T_p(f)=\sum_{m=0}^{p-1}\left({f(m)\over p}\right).$$
\begin{Lem}
\label{jacobsthal2} Let $p$ be an odd prime. Then
$$T_p(f)=\cases{-({a\over p}) &if $p\nmid ad$;\cr
p({c\over p}) &if $p|(a,d)$;\cr (p-1)({a\over p}) &otherwise.}
$$
\end{Lem}
{\it Proof}. If $p\nmid a$, then
$$({a\over p})T_p(f)=({4a\over
p})T_p(f)=\sum_{m=0}^{p-1}\left({(2am+b)^2-d\over p}\right)
=\sum_{m=0}^{p-1}\left({k^2-d\over p}\right),$$ where $k=2am+b$.
The proof is easily completed on invoking the previous lemma. (For
more details see, e.g., \cite[p. 79]{G}). \qed

\begin{Lem} Let $p$ be an odd prime. Then
$$a_p(f)=\cases{{-({a\over p})\over p-1-({d\over p})} &if $p\nmid ad$;\cr
0 &if $p|a,~p\nmid d$;\cr ({a\over p}) & if $p\nmid a,~p|d$;\cr
({c\over p}) & if $p|(a,d)$.}$$
\end{Lem}
{\it Proof}. The denominator in (\ref{definitie}) is easily
evaluated in prime arguments. On combining this computation with
Lemma \ref{jacobsthal2} the result
follows. \qed\\

\noindent The next result in the generic case was first
established by Andrew Granville (with a different proof).
\begin{Prop}
\label{lemmi} We have
$$a_D(f)=\cases{\displaystyle \left({c\over (D,a,d)}\right)
\left({a\over D/(D,a)}\right)\prod_{q|D\atop q\nmid ad} {-1\over
q-1-({d\over q})} &if $(D,a)|d$;\cr 0 &if $(D,a)\nmid d$.}$$
Alternatively,
$$a_D(f)=\left({c\over (D,a,d)}\right)\left({a\over
D/(D,a,d)}\right)\prod_{q|D\atop q\nmid ad} {-1\over q-1-({d\over
q})}.$$
\end{Prop}
{\it Proof}. Note that $a_D(f)=\prod_{p|D}a_p(f)$. Then invoke the
previous lemma. \qed\\

\noindent The latter result in combination with Theorem \ref{thm1}
gives:
\begin{Prop}
\label{gevolgen}
Assume Conjecture {\rm \ref{equidi}}. Let $f\in \cal F$.\\
{\rm 1)} If $\tau_D^-(f)\ne 0,1$, then $1/3\le \tau_D^-(f)\le 2/3$.\\
{\rm 2)} If $\tau_D^-(f)=0$ or $\tau_D^-(f)=1$, then $D|24ad$.
\end{Prop}
{\tt Remark} 1. We have $\tau_5^-(3X^2+7)=1/3$ and
$\tau_5^-(X^2+1)=2/3$ (thus the bounds in part 1 are sharp). One
computes that $\tau_{-3}^-(X^2+5)=1$ and thus $D|24ad$
in part 2 cannot be replaced by $D|8ad$.\\
{\tt Remark} 2. It can happen for a given $f\in \cal F$ that there
is no discriminant $D$ for which
$\tau_D^-(f)=1$, e.g. for $f(X)=X^2+X+41$.\\

\noindent The latter proposition strongly suggests that in order
to find large $c_g(f)$ we have to ensure that $\tau_D^-(f)=1$, where
$D$ denotes the discriminant of $\mathbb Q(\sqrt{g})$. This highly
restricts the possible choices of $D$. For Lehmer's polynomial
$L$, for example, one finds that $\tau_D^-(L)=1$ iff $D=-163,-3,24$
or $1304$.

\subsection{Higher degree $f$ }
If $f$ induces a permutation of $\mathbb F_p$ (that is, is
a permutation polynomial), then clearly
$a_p(f)=0$. E.g. if $f(X)=X^n+k$ and $(p-1,n)=1$, then $f$ induces
a permutation of $\mathbb F_p$ and hence $a_p(f)=0$.\\
\indent Suppose that $Y^2=f(X)$ is the Weierstrass equation of an elliptic
curve $E$ having conductor $N_E$. Hasse's
inequality yields $|a_p(f)|\le 2\sqrt{p}/(p-3)$ for $p>3$. 
It is well-known that $\sum_{j=0}^{p-1}({f(j)\over p})$
is the trace of Frobenius over $\mathbb F_p$. In the remainder
of this section it is assumed that the conditions of Theorem \ref{thm1} are satisfied, so
that Theorem \ref{thm1} can be invoked.
It follows that if $D\equiv 1({\rm mod~}4)$ and $(N_E,D)=1$, then $\tau_D^{-}(f)=1/2$
iff there is prime $p$ dividing $D$ such that $E$ is supersingular at $p$. Since
Deuring it is known that the number of supersingular primes $p\le x$ in case of
a CM curve E grows asymptotically as $\pi(x)/2$ and hence for almost all quadratic
fields of odd discriminant $D$ one has in this case $\tau_D^{-1}(f)=1/2$ (again under the
conditions of Theorem \ref{thm1}). On the other hand, if $E$ does not have complex
multiplication one finds using the result of Serre that the number of supersingular primes $p\le x$
is then bounded by $\ll x(\log x)^{-5/4+\epsilon}$ that for a positive proportion of the fundamental
discriminants $D\equiv 1({\rm mod~}4)$ one has $\tau_D^-(f)=1/2$. 

\section{Heuristics for the proportion of primitive roots}
In the previous section we gave an heuristic for the
proportion $\tau_D^-(f)$ of
primes $p=f(n)$ such that $({g\over p})=-1$. In this section we do
the same but with the more stringent condition that $g$ should be
a primitive root modulo $p$. Numerical work suggests the truth of:
\begin{Conjecture}
Suppose that $f(X)\in \mathbb Z[X]$ represents infinitely many
primes. Then the quotient of
$$\# \{p\le x:f(m)=p{\rm ~for~some~}m{\rm ~and~}g{\rm
~is~a~primitive~root~mod~}p\}$$ and $\# \{p\le x:f(m)=p{\rm
~for~some~}m\}$ tends to a limit as $x$ tends to infinity, that is
the relative proportion of primes $p$ such that $g$ is a primitive
root mod $p$ and moreover $p$ is represented by $f(x)$ exists. Let
us denote this conjectural density by $\delta_g(f)$.
\end{Conjecture}
\indent In the remainder of this section it is assumed that the
latter conjecture holds true. It is also supposed that $g$ is not
an $h$th power of an integer for any $h\ge 2$.\\
\indent Suppose that $g$ is such that $\tau_D^-(f)=1$, where $D$ is
the discriminant of $\mathbb Q(\sqrt{g})$ (the most relevant
case for our purposes). Then, by an argument similar to that
used in the derivation of (\ref{p1}), one is led to believe that a
good approximation for $\delta_g(f)$ should
be
\begin{equation}
\label{deltaf}
\delta(f):=\prod_{q>2}\left(1-{\# \{s({\rm mod~}q):f(s)\equiv
1({\rm mod~}q)\} \over q\#\{s({\rm mod~}q):f(s)\not\equiv 0({\rm
mod~}q)\}}\right).
\end{equation}
In case $f(X)=AX^2+B$
a short calculation shows that
$$\delta(f)=\prod_{q|(A,B-1)\atop q>2}(1-{1\over q})
\prod_{q\nmid 2A}\left(1-{\{1+\left({-A(B-1)\over q}\right)\}\over
q(q-1-\left({-AB\over q}\right))}\right).$$
For general quadratic $f(X)=aX^2+bX+c$ one finds that
\begin{equation}
\label{afschatting2}
{\varphi((a,b,c-1))\over (a,b,c-1)L(2,(d/.))} \ll \delta(f)\ll {\varphi((a,b,c-1))\over (a,b,c-1)L(2,(d/.))},
\end{equation}
where $d=b^2-4a(c-1)$.\\
\indent If $\delta(f)$ is close to $1$, then
$$\delta_1(f):=\prod_{q|(A,B-1)\atop q>2}(1-{1\over q})
\prod_{q\nmid 2A}\left(1-{\{1+\left({-A(B-1)\over q}\right)\}\over
q^2}\right),$$
yields a quite good approximation to $\delta(f)$; compare (\ref{lehmertje})
with (\ref{p1}).
Clearly the idea in finding a large value of $c_g(f)$ is to find
$f$ such that $\delta(f)$ is close to $1$. For this results from
the theory of prime producing quadratics can be used.

\section{Prime producing quadratics}
Let $f_A(X)=X^2+X+A$, with $A>0$ a positive integer. Euler
discovered in 1772 that $X^2+X+41$ satisfies
$\pi_{f_{41}}(39)=40$. It can be shown that $\pi_{f_A}(A-2)=A-1$
iff $A\in \{2,3,5,11,17,41\}$, see Mollin \cite{Mollin}, and that
this is related to the class number one problem. The connection with the
class number one problem dates back to Frobenius (1912) and Rabinowitsch (1913). The discriminant
of $f_A(X)$ is given by $\Delta=1-4A$. Note that if $A$ is even,
then $2|f_A(x)$ and so we may assume that $A$ is odd and hence
$\Delta\equiv 5({\rm mod~}8)$. If for a prime $q$, $({\Delta\over
q})=-1$, then the values of $f_A$ are not divisible by $q$. So if
$({\Delta\over q})=-1$ for many consecutive primes $q$, the values
of $f_A$ have a better chance of being prime, in particular if
$\Delta$ is also small. Thus we want
\begin{equation}
\label{L} L(1,\chi)=\prod_q {1\over 1-\chi(q)/q},
\end{equation}
where $\chi_{\Delta}(n)=(\Delta/n)$ and $(./n)$ is the Kronecker symbol to
be small. Since with two exceptions $\pi
h/\sqrt{|\Delta|}=L(1,\chi_{\Delta})$, we want the class number $h$ to be
small. By (\ref{HLL}) one should have, as $x$ tends to infinity,
$\pi_{f_A}(x)\sim C(\Delta)x/\log x$, where
$$C(\Delta)=\prod_{q\ge 3}\left(1-{({\Delta\over  q})\over
q-1}\right).$$ It is easy to show (using that $(\Delta/2)=-1$)
that
\begin{equation}
\label{nietvanmij} C(\Delta)={\zeta(4)\over 2
L(1,\chi_{\Delta})L(2,\chi_{\Delta})}\prod_{q|\Delta}(1-{1\over q^4}) \prod_{q\ge
3\atop ({\Delta\over q})=1} \left(1-{2\over q(q-1)^2}\right).
\end{equation}
Shanks has computed $C(-163)=3.3197732 \ldots$ and
$C(-111763)=3.6319998 \ldots$. Thus Beeger's \cite{B} polynomial
$X^2+X+27941$ should produce asymptotically more primes than
Euler's. One computes that $\pi_{f_{41}}(10^6)=261080$ and
$\pi_{f_{27941}}(10^6)=286128$. On the other hand $\pi_{f_{41}}(39)=40$,
whereas $\pi_{f_{27941}}(39)=30$. The constant $C(\Delta)$ can become
arbitrarily large: for every $\epsilon>0$ there are infinitely
many $\Delta$ such that
$$(1/2+\epsilon)e^{\gamma}\log \log |\Delta| < C(\Delta)
< (1+\epsilon)e^{\gamma}\log \log |\Delta|,$$
where $\gamma$ denotes Euler's constant (see \cite[p. 511-512]{JW}).\\
\indent Quadratics that produce
too many primes contradict the Generalized Riemann Hypothesis. If there are
lots of Siegel zeros this can be used to infer results on the growth of $\pi_f(x)$. This is
akin to Heath-Brown's result that if there are many Siegel zeros, then the twin primes behave
as expected. For more on the analytic aspects of prime-producing polynomials, see \cite{GM}.\\
\indent In order to find $\Delta$ with $({\Delta\over q})=-1$ for
many consecutive primes $q$, special purpose devices have been
built (some even involving bicycle chains !). For a nice account
of this see Lukes, Patterson and
Williams \cite{LPW}.\\
\indent In searching for good prime producing quadratics it is
thus tantamount to find $\Delta$ for which $C(\Delta)$ is large.
Similarly, for Problem 1 we want $\delta(f)$ to be close to
$1$. Equation (\ref{nietvanmij}) shows that finding a large value of $C(\Delta)$ 
amounts to finding $\Delta$ such that $L(1,\chi_{\Delta})$ is small. For our problem
at hand, however, the issue is rather to find small $L(2,\chi_{\Delta})$. 
To see this note that $\delta_1(f)$ is a rational multiple of
\begin{equation}
\prod_{q\ge 3}\left(1-{\{1+({\Delta\over q})\}\over q^2}\right)={3\over 4}\zeta(2)\prod_{q\ge 3}\left(1-{({\Delta\over q})\over 
q^2-1}\right).
\end{equation}
It is not difficult to show that for Re$(s)\ge 1$ 
\begin{equation}
\label{mijexpressie}
\prod_{q\ge 3}\left(1-{\chi_{\Delta}(q)\over q^s-1}\right)=\epsilon(s)
{\zeta(2s)\over L(s,\chi_{\Delta})}\prod_{q|\Delta}(1-{1\over q^{2s}})\prod_{q\ge 3\atop ({\Delta\over q})=1}
\left(1-{2\over q^s(q^s-1)}\right),
\end{equation}
where $\epsilon(s)=1+2^{-s}({\Delta\over 2})$. For $s=1$ we obtain an expression for $C(\Delta)$ and
for $s=2$ we obtain an expression closely related to $\delta_1(f)$.
In case $s=1$ the latter product in the expression does not converge very well and
preference is to be given to expression (\ref{nietvanmij}). However, in case $s=2$ 
expression (\ref{mijexpressie}) is quite usable. The special value $L(2,\chi_{\Delta})$ involved
can be evaluated with high precision, see \cite{JW}.\\
\indent Let $\alpha\ge 1$. If $f(X)$ is a prime producing quadratic, then $g_{\alpha}(X)=2^{\alpha}f(X)+1$
is likely to be primitive root producing for those $g$ satisfying $\tau_D^{-}(g_{\alpha})=1$, with $D$ the
discriminant of $\mathbb Q(\sqrt{g})$. Conversely, if $g(X)$ is a primitive root producing quadratic, 
then we can write $g(X)-1=2^{\alpha}(aX^2+bX+c)$ with $\alpha\ge 0$ and $(a,b,c)=1$. Write
$h(X)=aX^2+bX+c$.
If $N_2(h)=0$, then $h$ is likely to be prime producing. Thus the connection
between primitive root producing and prime producing quadratics is rather intimate.

\section{Finding primitive root producing quadratics}
In general an approach to Problem \ref{P1} is to find a small
integer $d$ such that $({d\over q})\ne 1$ for as many small odd
primes $q$ as possible. Thus we hope to ensure that $\delta(f)$
(the quality of $f$) is very close to $1$. We factorize $d$ as
$d_1d_2$ and choose a small $\alpha$. Then we consider primes $p$
of the form $2^{\alpha}d_1n^2+2^{\alpha}d_2+1$. Since we want
$({g\over p})\ne 1$ for all primes of the latter form, the choice
of $g$ is rather restricted: under Conjecture \ref{equidi} the
discriminant $\mathbb Q(\sqrt{g})$ has to be a divisor of $24
d_1(2^{\alpha}d_2+1)$
by Proposition \ref{gevolgen}.
It can happen that no suitable $g$ can be found and then $\alpha$
can be adjusted. If $g$ has the required property, so has $k^2g$
for every integer $k$. Now we vary over $k$ in the hope of finding
a large value of $c_{k^2g}(2^{\alpha}d_1X^2+2^{\alpha}d_2+1)$.
Another variation option we have is to consider primes $p$ of the
form $2^{\alpha}d_1r_1n^2+2^{\alpha}d_2r_2+1$ with $r_1r_2$ a
square and with $r_1r_2$ having only large prime factors. The
corresponding value of $\delta(f)$ changes little by this and again
we can search for a large value of
$c_g(2^{\alpha}d_1r_1X^2+2^{\alpha}d_2r_2+1)$. 
(In this variation $g$ remains fixed and thus it can be used in 
dealing with Problem 2.)
Since we want
$({g\over p})\ne 1$ usually some mild congruence conditions on
$r_1$ and $r_2$ have to be imposed. A further variation
possibility is to replace $n$ by $\gamma n+\delta$. However,
computational practice suggests this is only effective when $\gamma=1$.\\
\indent The asymptotic (\ref{oudprobleempje}) suggests that it is
crucial to get a large value of $\delta(f)$: if this value is
not close enough to $1$, then there is not much to be gained by
letting $k$ run over a large range (note that in general
$p_1=\delta_g(f)$).\\

\noindent {\tt Example 1}. The number $d=4472988326827347533$
satisfies $(d/p)=-1$ for the primes $p=3,\dots,283$ by Table 4.3
of \cite{JW}. A factor of $d$ is $d_1=252017$. Let $d_2=d/d_1$.
Let $f(X)=1008068X^2+16921429448X+15753313937$. (This is
$4d_1(X+8393)^2-4d_2+1$.)  The first `bad' prime equals
$432050978399143373$. It turns out that $c_{170363492}(f)=22779$.
One finds that $\delta(f)\approx 0.999453$
and that $M(\delta(f),145700)\approx 22779$.\\

\noindent {\tt Example 2}. (Y. Gallot). We let $d$ be as in Example 1, 
$d_1=230849$ and $d_2=d/d_1$. Let
$f(X)=64d_1(X+728069)^2-64d_2+1$ and $g=17^{2}\cdot 230849=66715361$.
Then $c_g(f)=25581$. This is the presently largest known value of $c_g(f)$ for an $f$
having positive discriminant. One finds that $\delta(f)\approx 0.999453$
and that $M(\delta(f),675200)\approx 25581$.\\
\indent Let $f(X)=64d_1(X+56943)^2-64d_2+1$. Then $d_{24}(f)=21690$. This
is the record for $c_g(f)$ with $|g|<100$ (cf. Problem 2).\\

\noindent {\tt Example 3}. The number $d=9828323860172600203$
satisfies $(-d/p)=-1$ for the primes $p=3,\dots,277$ by Table 4.1
of \cite{JW}. A factor of $d$ is $d_1=54151$. Let $d_2=d/d_1$.
Let $f(X)=866416X^2+2903975582404049$. (This is
$16d_1X^2+16d_2+1$.)  It turns out that $c_{23731350844}(f)=18176$.
Let $f_1(X)=f(X+599206)$. One computes that $c_{72922}(f_1)=29083$.
Let $f_2(X)=d_1(X+1484224)^2+d_2+1$. Then $c_{17431902}(f_2)=31082$.
This is the presently largest known value of $c_g(f)$ for an $f$
having negative discriminant and was discovered by Yves Gallot.
One finds that $\delta(f_2)\approx 0.999535$
and that $M(\delta(f_2),1066000)\approx 31082$.

\section{On the (un)boundedness of $c_g(f)$}
A tool in investigating this is an extension of a criterion of
Chebyshev which is discussed in the next section.

\subsection{Extension of a primitive root criterion of\\ Chebyshev}
It is an old result of Chebyshev that if $p_1\equiv 1({\rm
mod~}4)$ is prime and $p_2=2p_1+1$ is also prime, then $g=2$ is a
primitive root modulo $p_2$. Under TC(2) it then follows that $2$
is a primitive  root for  infinitely many primes. Already in the
19th century Chebyshev's criterion was extended to some numbers
other than $2$, see e.g. \cite{W}. In this section an analogue of
Chebyshev's criterion is derived for every integer $g\in G$. This
criterion plays a
keyrole in the proof of Theorem \ref{AG}.\\
\indent It is not known whether there are infinitely many primes
satisfying Chebyshev's criterion, but it can be shown that there
are infinitely many primes satisfying a somewhat weaker version of
it. This can then be used to show, e.g., that at least  one of the
numbers $2,3$ and $5$ is a primitive root for infinitely many
primes \cite{HB}.

\begin{Lem}
\label{previous} Let $g\ge 3$ be an odd squarefree integer. There
exists an integer $a$ such that $(a,g)=1$ and $({8a+1\over
g})=-1$.
\end{Lem}
{\it Proof}. It is easy to see that the result holds true in case
$g$ is an odd prime. In case $g\ge 5$ is an odd prime, likewise
there exists an integer $b$ such that $(b,g)=1$ and
$({8b+1\over g})=1$. From these two observations the result
follows
on invoking the Chinese Remainder Theorem. \qed

\begin{Lem}
\label{chebbie} Suppose that $g\in G$. Write $g={g_0}^2g_1$ with
$g_1$ squarefree.
Let $g_2=|g_1|$ if $g_1$ is odd and $g_2=|g_1/2|$ otherwise.\\
\indent For parts {\rm 1} and {\rm 2} it is assumed that $g_1\ne \pm 2$.\\
{\rm 1)} Let $a$ be any integer such that $(a,g_2)=1$ and
$({8a+1\over g_2})=-1$ (by Lemma {\rm \ref{previous}} at least one
such integer exists). If $p_1$ is a prime of the form $g_2k+a$
such that $p_2:=8p_1+1$ is also a prime and $g^8\not\equiv
0,1({\rm mod~}p_2)$, then
$g$ is a primitive root modulo $p_2$.\\
{\rm 2)} Under {\rm TC(2)} there are infinitely many primes $p_1$
satisfying the
conditions of part {\rm 1}.\\
{\rm 3)} Assume that $g_1=\pm 2$. If $p_1$ is a prime and
$p_2:=2p_1+1$ is a prime, then $g$ is a primitive root modulo
$p_2$ if $p_1\equiv {\rm sgn}(g)({\rm mod~}4)$ and $g^2\not \equiv
0,1({\rm mod~}p_2)$. If {\rm TC(2)} holds true, there are
infinitely many primes $p$ such that $g$ is a primitive root
modulo $p$.
\end{Lem}
{\it Proof}. 1) The assumption $g^8\not\equiv 0,1({\rm mod~}p_2)$
ensures that the order of $g$ modulo $p_2$ exists and is a
multiple of $p_2$. Since
$$({g\over p_2})=({g_1\over p_2})=({g_2\over p_2})=({p_2\over
g_2})=({8a+1\over g_2})=-1,$$ and $-1=({g\over p_2})\equiv
g^{4p_1}({\rm mod~}p_2)$, the order must be
$8p_1=p_2-1$.\\
2) We have to show that for each prime $p$ there exists $k$ for
which
\begin{equation}
\label{puh} (g_2k+a)(8g_2k+8a+1)\not\equiv 0({\rm mod~}p).
\end{equation}
For $p=2$ this is clear. In case $p|g_2$ this follows since we have
$(a,g_2)=1$ and $(8a+1,g_2)=1$. For the remaining primes $p$
there are at least $p-2\ge 1$ choices  of $0\le k<p$ such that
(\ref{puh}) is
satisfied.\\
3) Similar to the proof of parts 2 and 3.\qed

\begin{Cor}
Artin's primitive root conjecture is true, assuming {\rm TC(2)}.
\end{Cor}
Recall that Artin's conjecture (1927) asserts that any integer
$g\in G$ is a primitive
root for infinitely many primes $p$.\\
\indent Another generalisation of Chebyshev's criterion is in the
direction of cubic reciprocity. For example, if $p$ is an
odd prime such that $q=1+6p$ is a prime then $3$ is not a
primitive root mod $q$ iff we can write $4p=n^2+243m^2$ with $n,m$
integers. This criterion is due to Fueter \cite{F}.

\subsection{A conditional result on $c_g(f)$}
 \noindent Lemma \ref{chebbie} will be used in the proof of the
 following theorem, the basic idea of which is due to Andrew Granville.
\begin{Thm}
\label{AG} Let $N\ge 1$ be an integer. Assume {\rm TC($2N$)}.
Suppose that $g\in G$. Then there exist integers $A_1$ and $C_1$
such that $A_1n^2+C_1$ is prime for $n=1,\dots,N$ and $g$ is a
primitive root for each of these primes.
\end{Thm}
Here and in the sequel $A_1$ and $C_1$ are allowed to depend on
$N$.
\begin{Cor} Assume {\rm TC($2N$)} for every $N\ge 1$.
Let $g\in G$ be fixed. The number $c_g(AX^2+C)$ can be larger than any
prescribed number.
\end{Cor}
{\tt Remark}. Let $N\ge 1$ be an integer and $g\in G$. Perhaps it is possible to show under TC that there exist 
integers $A_1$ and $C_1$ such that $A_1n^2+C_1$ is prime for $n=1,2,\ldots,N+1$ and $g$ is a primitive
root for the first $N$ of these primes, but not for the $N+1$th. This would show that $c_g(AX^2+C)$ can
assume any prescribed natural number as value under TC.\\

\noindent {\it Proof of Theorem} \ref{AG}. We adopt the notation of Lemma \ref{chebbie}
and assume that $g_1=\pm 2$ (the remaining case being similar).\\ 
\indent Let $A = \prod_{p\le
2N}p$ and $C$ be the smallest integer $>2N$ with $C\equiv a({\rm
mod~}g_2)$ for which $C$ and $8C+1$ are both primes ($C$ exists by
part 2 of Lemma \ref{chebbie}). Consider the $2N$-tuplet of
numbers $g_2At+C+g_2An^2$ for $n=1,\dots,N$ and
$8g_2At+8C+1+8g_2An^2$ for $n=1,\dots,N$ for integer $t$. TC($2N$)
predicts that there will be infinitely many $t$ for which these
are all prime, provided there is no obstruction modulo a prime $p$
(i.e. it is not true that for every $t$ at least one of the forms
is divisble by $p$). (We will take $A_1=8g_2A$ and
$C_1=8g_2At+8C+1$ above for one of these $t$'s such that,
moreover, none of the primes $p(n)$ of the form $A_1n^2+C_1$ with
$n=1,\dots,N$ satisfies $g^8\equiv 0,1({\rm mod~}p(n))$). Now for
$p\le 2N$, we see that $p|A$ and $p\nmid C(8C+1)$, so $p$ never
divides any of the forms. If $p|g_2$ the first $N$ forms are
$\equiv a({\rm mod~}p)$ and the second $N$ forms are $\equiv
8a+1({\rm mod~}p)$. The conditions on $a$ ensure that
$a(8a+1)\not\equiv 0({\rm mod~}p)$. In general there are at most
$2N$ values  of $t$ for which at least one of our $2N$ linear
forms is divisible by $p$, so if $p>2N$ and $p\nmid g_1$, there
exists an
integer $t$ such that none of them is divisible by $p$.\\
\indent Let $p(n)=A_1n^2+C_1$. Now for $1\le n\le N$ each $p(n)$
is a prime for which $(p(n)-1)/8$ is also a prime and satisfies
the conditions of part 1 of Lemma \ref{chebbie} and hence
$g$ is a primitive root modulo $p(n)$. \qed

\begin{Lem}
\label{vorige} Suppose that $g_i\ne -1$ for $i=1,\dots,s$ and that
\begin{equation}
\label{minus} ({g_1\over p})=\dots=({g_s\over p})=-1
\end{equation}
for infinitely many primes $p\equiv 2({\rm mod~}3)$, then there
exists $1\le m\le 2$, $a$ and $f$ with $(a,f)=1$, such that for
every prime $q$ satisfying $q\equiv a({\rm mod~}f)$  for which
$q_1=2^mq+1$ is also a prime and $g_i^{2^m}\not\equiv 0,1({\rm
mod~}q_1)$ for $i=1,\dots,s$, then the integers $g_1,\dots,g_s$ are
simultaneously
primitive roots modulo $q_1$.
\end{Lem}
{\it Proof}. Let $Q=\{q_1,\dots,q_t\}$ be the set of odd primes
dividing the discriminant of $\mathbb Q(\sqrt{g_i})$ for some
$1\le i\le s$. Let $A_{+1}(q)$ be the set of non-zero quadratic
residues modulo $q$ and $A_{-1}(q)$ the set of quadratic
non-residues. It is a consequence of quadratic reciprocity that
there exist $\epsilon_i\in \{-1,1\}$ with the property that for
each choice of elements $\alpha(\epsilon_i)\in A_{\epsilon_i}(q)$,
there are infinitely many primes $p$ satisfying (\ref{minus}) such
that, moreover, $p\equiv \alpha(\epsilon_i)({\rm mod~}q_i)$ for
$1\le i\le t$. The condition that $p\equiv 2({\rm mod~}3)$ now
ensures that we can pick $\alpha(\epsilon_i)\ne 1$. The argument
can easily be extended to take the behaviour at the prime two into
account. One sees one can pick $\beta\in \{3,5,7\}$ such that
there are infinitely many primes $p$ satisfying (\ref{minus}) such
that $p\equiv \beta({\rm mod~}8)$ and $p\equiv
\alpha(\epsilon_i)({\rm mod~}q_i)$ for $1\le i\le t$. Setting
$f=8q_1\cdots q_t$, one then finds that $a$ with $2^ma+1\equiv
\beta ({\rm mod~}8)$ and $2^ma+1\equiv \alpha(\epsilon_i) ({\rm
mod~}q_i)$ for $1\le i\le t$ exists and satisfies the requirement
$(a,f)=1$, provided we set $m=2$ if $\beta=5$ and $m=1$
otherwise. The proof is then finished by an
argument as used in the proof of Lemma \ref{chebbie}. \qed\\

\noindent The following result generalizes Theorem \ref{AG}.
\begin{Thm}
Let $s\ge 1$ be an integer and let $g_1,\dots,g_s$ be integers
$\ne -1,0,1$. Let $0\le e_1,\dots,e_s\le 1$. Suppose that
$\prod_{i=1}^s g_i^{e_i}$ is not a square if $e_1+\ldots+e_s$ is
odd. Suppose furthermore that the discriminant of each of the
fields $\mathbb Q(\sqrt{g_i})$ is not divisible by $3$. Then there
exist integers $A$ and $C$ such that $p(j)=Aj^2+C$ is prime for
$1\le j\le n$ and each of the $g_i$ is a primitive root modulo
$p(j)$.
\end{Thm}
{\it Proof}. Using the argument at p. 37 of Heath-Brown \cite{HB},
one easily infers that the conditions of Lemma \ref{vorige} are
satisfied. Thus there exist numbers $a,f$ and $m$ as in that
lemma. Now proceed as in the proof of Theorem \ref{AG}. Thus take
$C$ to be the smallest integer $>2N$ with $C\equiv a({\rm mod~}f)$
and replace $8C+1$ by $2^mC+1$. The rest of
the argument is left as a (copy) exercise to the interested reader. \qed\\

\noindent {\tt Remark}. I do not see how to prove this result with
for example $g_1=-25$ and $g_2=3$, although in this case under GRH
it can be shown that there are infinitely many primes $p$ such
that both are primitive roots \cite{M}. In essence the question amounts to
this one: for each $N\ge 1$ are there $A$ and $C$ such that
$p(j)=Aj^2+C\equiv 7({\rm mod~}12)$ are all prime and $3$ is a
primitive root mod $p(j)$ for $1\le j\le N$ ? One seems to be
forced to use cubic reciprocity, cf. Fueter's criterion (Section 9.1).

\section{Conclusion}
By {\it Griffin's dream} I understand the dream to find a
polynomial $f$ that represents infinitely many distinct primes and an
integer $g$ such that for all primes $p=f(n)$ with $p\nmid g$ and $n\ge
0$, the integer
$g$ is a primitive root modulo $p$.
\begin{Conjecture} \label{weereen}$~$\\
{\rm 1)} For quadratic $f$ Griffin's dream cannot be realized, i.e.
$c_g(f)<\infty$.\\
{\rm 2)} Let $m\ge 1$ be arbitrary. For $g\in G$ there exist $f$ such
that $c_g(f)>m$.
\end{Conjecture}
I base part 1 on the following proposition and the observation
that if an event can
 occur with positive probability it will eventually occur (after enough
repetition).
\begin{Prop}
Let $f\in \mathbb Z[X]$ be quadratic. Then $\delta(f)<1$.
\end{Prop}
{\it Proof}. Suppose that $\delta(f)=1$. Then from (\ref{deltaf}) one infers the existence of a
fundamental discriminant $\Delta$ such that $({\Delta\over q})=-1$
for all but finitely many primes $q$. Since $\prod_{p\le x}(1+1/p)\sim c_1\log x$ for some
$c_1>0$ by a result of Mertens, it then follows from (\ref{L}) that
$L(1,\chi_{\Delta})=0$. However,
$L(1,\chi_{\Delta})>0$ as is well-known. \qed\\

\noindent The motivation for part 2 of Conjecture \ref{weereen} is provided by Theorem \ref{AG}.\\
\indent Whereas the problem of finding prime producing polynomials amounts to  finding 
$D$ for which $L(1,\chi_{D})$ is small (cf. the estimate (\ref{afschatting1})), the problem of finding 
primitive root producing polynomials amounts to finding
$D$ for which $L(2,\chi_{D})$ is small (cf. the estimate (\ref{afschatting2})).\\

\noindent {\bf Acknowledgement}. I'd like to thank Bruce Berndt
for pointing out reference \cite{G} in relation with Lemma
\ref{Jacobsthal}. Igor Spharlinski kindly sketched a proof of the
latter lemma. As so  often Yves Gallot and Paul Tegelaar kindly provided computational assistance.\\ 
\indent Special thanks are due to Andrew Granville, if it
were not for him the paper would have looked quite differently:
especially Section 9 would not have been there. Thanks are also
due to Richard Mollin for passing on a question of mine to Andrew.

\end{document}